\documentclass[12pt,leqno]{article}
\usepackage{amsmath,amsfonts}
\newtheorem{lem}{Lemma}[section]

\newtheorem{theo}[lem]{Theorem}
\newtheorem{cor}[lem]{Corollary}
\newtheorem{prop}[lem]{Proposition}
\newcommand{\vers}{\mathop{\longrightarrow}} 
 
\newcommand{\R}{{\mathbb R}}

\numberwithin{equation}{section} 

\title{Some explicit Krein representations of certain subordinators, including the Gamma process}
\author{C. Donati-Martin and M. Yor \thanks
{Laboratoire de Probabilit\'es et Mod\`eles Al\'eatoires, Universit\'e Paris 6, Site Chevaleret, 13 rue Clisson, F-75013 Paris. email: donati@ccr.jussieu.fr }}
\date{}
\begin{document}
\maketitle 
\begin{abstract}
We give a representation of the Gamma subordinator as a Krein functional of Brownian motion, using the known representations for stable subordinators and  Esscher transforms. In particular, we have obtained Krein representations of the subordinators which govern the two parameter Poisson-Dirichlet 
family of distributions \cite{PY97}.
\end{abstract}
\section{Introduction} 
{\bf (1.a)} The aim of this paper is to represent explicitly a particular class of subordinators $(S_t; t \geq 0)$, i.e. increasing $\R_+$-valued L\'evy processes, as  inverse local times of $\R_+$-valued diffusions. This problem was raised in It\^o-Mc Kean \cite{IMK}, and completely solved, in a theoretical manner, thanks to Krein's representation of strings, by Knight \cite{K} and  Kotani-Watanabe \cite{KW} independently and simultaneously in 1981-1982.

\noindent
Roughly, if the L\'evy measure $\nu$ of $(S_t; t \geq 0)$ admits a density $h$ with respect to the Lebesgue measure: $\nu(dy) = h(y) dy$ and if $h(y) = \int_0^\infty \mu(dx) \exp(-yx)$ for some positive measure $\mu$, then the above mentioned problem, which we shall call the Krein representation problem, is solved in the affirmative. See, besides  \cite{K} and \cite{KW}, a number of other papers dealing with this question, e.g. Bertoin (\cite{B}, \cite{B1}), K\"uchler \cite{Ku}, K\"uchler and Salminen \cite{KuS}.

\vspace{.3cm}
\noindent
{\bf (1.b)} In this paper, rather than discussing this problem in general, we solve it for the 
 two parameter family of subordinators $(S_t^{\alpha, \beta}; t \geq 0)$  whose L\'evy measures are 
$$\nu_{\alpha, \beta}(dy) =  C \frac{ \exp(- \beta y)}{y^{\alpha +1}} dy, \ 0\leq \alpha <1, \beta \geq 0$$
(for $\alpha = 0$, $\beta$ is assumed $>0$). \\
The constant $C$ plays a simple role and might, a priori, be suppressed from our discussion; however, it is in fact very helpful to keep  this further degree of freedom in order to consider the most convenient  local time at 0 (for the underlying diffusion), which, as is well known, may be chosen up to a multiplicative constant.

\vspace{.3cm}
\noindent
{\bf (1.c)} It goes back at least to Molchanov-Ostrovski \cite{MO} that for $0< \alpha<1$,  $S_t^{\alpha,0}$, the  stable subordinator of index $\alpha$, may be realized as 
 the inverse local time of a Bessel process (which we shall denote as $BES(-\alpha)$) of dimension $\delta = 2(1- \alpha)$.
 Since, by inspection of the L\'evy measures: $\nu_{\alpha, \beta}(dy) =    \exp(- \beta y) \nu_{\alpha, 0}(dy)$,  $S_t^{\alpha, \beta}$  may be obtained as an Esscher transform of $S_t^{\alpha, 0}$, it is natural to look for a Girsanov transform of $BES(-\alpha)$ whose inverse local time is distributed as
 $S_t^{\alpha, \beta}$. As we show below, this is indeed the case with $BES(-\alpha, \beta\!\downarrow)$ the downwards $BES(-\alpha)$ process with "drift" $\beta$, following the terminology of Watanabe \cite{W1}, \cite{W2} and Pitman-Yor \cite{PY1}. \\
  Finally, the case $\alpha = 0$, $\beta >0$, which corresponds to the Gamma process 
$(S_t^{0, \beta}; t \geq 0)$ is obtained from the case $\alpha >0$ by letting $\alpha \vers 0$ in a suitable manner.

\vspace{.3cm}
\noindent
{\bf (1.d)} It is also natural to look for some representation of these subordinators as continuous additive functionals of a Brownian motion $(B_s; s \geq 0)$ taken at the inverse local time $(\tau^*_t; t \geq 0)$ of that Brownian motion. This is done expressing  the Bessel processes with drift in terms of their Feller representations, using scale functions and speed measures. For example, it was remarked in Biane-Yor \cite{BY} that:
\begin{equation} \label{repstable}
( \int_0^{\tau^*_t} |B_u|^{\frac{1}{\alpha}-2} \ du, \; t \geq 0 ) \stackrel{(law)}{=} (S_t^{\alpha,0}; t \geq 0),
\end{equation}
and that all symmetric stable L\'evy processes may be obtained in this manner, replacing the even power
$|x|^{\frac{1}{\alpha}-2}$  by "symmetric" powers i.e.: $\sigma(x) = sgn(x) |x|^{\frac{1}{\alpha}-2}$. For example,
$(\frac{1}{\pi} \int_0^{\tau^*_t} \frac{ds}{B_s}; t \geq0)$ is a standard Cauchy process. More generally, every asymmetric stable L\'evy process may be represented in a similar  manner from Brownian motion. Although in the sequel, we shall also extend \eqref{repstable} to present in general   $(S_t^{\alpha, \beta}; t \geq 0)$  in terms of Brownian additive functionals taken at $(\tau^*_t; t \geq 0)$, these presentations are not so simple, and we prefer to those the Krein representations evoked in {\bf (1.c)}.

\vspace{.3cm}
\noindent
{\bf (1.e)} The rest of the paper is organized as follows. Our main Krein representation results are presented in Section 2. The full proofs are given in Section 3. The Brownian additive functional representations  are discussed in Section 4. Finally, in Section 5, we also give some Krein representations of the symmetric L\'evy processes on $\R$ (without Gaussian component) whose 
L\'evy measures are given by:
$$\tilde{\nu}_{\alpha, \beta} (dy) = \frac{\exp(-\beta |y|) }{ |y|^{\alpha +1}} dy$$
($ 0 \leq \alpha<2, \ \beta >0$).

\vspace{.3cm}
\noindent
{\bf (1.f)} As an end to this introduction, let us point out that, the subordinators $S^{\alpha, \beta}$ we represent here, and their symmetric counterparts, are arguably the most studied and used among L\'evy processes, and this, for the following reasons: for $\alpha >0$, $S^{\alpha, \beta}$ is obtained by Esscher transform (see \cite{Es, Sh}) from the fundamental stable $(\alpha)$ subordinator, hence it "retains" some scaling property, while the Gamma process (\cite{EY}, \cite{VY}, \cite{TV}, \cite{TVY}, \cite{TVY1}) and the variance-gamma processes (\cite{MS}, \cite{MM}, \cite{MCC}) have some fundamental quasi-invariance properties, which make them comparable, in some respect, to Brownian motion with drift.

\vspace{.3cm}
\noindent
{\bf (1.g)} Finally, we refer the reader to a systematic compendium \cite{DRVY} of constants $C_\alpha$ related to various choices found in the literature of local times with respect to $BES(-\alpha)$. We also
intend in \cite{DY} to study a number of properties of the symmetrized $BES(\alpha, \beta\!\downarrow)$ processes with a view towards a discussion of Krein representations for the variance-gamma processes.
%%%%%%%%%%%%%%
%%%%%%%%%%%%%%%
\section{Definition of the $BES(-\alpha, \beta\!\downarrow)$ processes and main Krein representation results}
{\bf (2.1)} We first recall that for $0< \alpha<1$, $BES(-\alpha)$, which we call the Bessel process of index $-\alpha$, (or dimension $d= 2(1-\alpha)$) is the $\R_+$-valued
 diffusion with infinitesimal generator
$$\frac{1}{2} \frac{d^2}{dx^2} + \frac{1-2 \alpha}{2x} \frac{d}{dx},$$ which is instantaneously reflecting at 0. We denote by $P_x^{(-\alpha)}$ its distribution on $C(\R_+, \R_+)$, where $R_t(\omega) = \omega(t)$,  ${\cal R}_t = \sigma\{R_s; s \leq t\}$ and $R_0 = x$. \\
For the sequel, it  will be  convenient to introduce the parameter $\theta = \sqrt{2\beta}$ for $\beta >0$. We may now define 
$BES(-\alpha, \beta\!\downarrow)$ as the diffusion  with law $(P_x^{(-\alpha), \beta\!\downarrow})$, obtained by Girsanov\footnote{We shall also use the well-known general fact that such an absolute continuity relationship extends with $t$ replaced by any stopping time $T$ on the set $(T <\infty)$.} transform from $BES(-\alpha)$:
\begin{equation} \label{defbesdrift}
P_x^{(-\alpha), \beta\!\downarrow}|_{{\cal R}_t} = \frac{\hat{K}_\alpha (\theta R_t)}{\hat{K}_\alpha (\theta x)} \exp( C_\alpha (2\beta)^\alpha l_t - \beta t) P_x^{(-\alpha)}|_{{\cal R}_t}
\end{equation}
where $\hat{K}_\alpha (x) = x^\alpha K_\alpha (x), x \geq 0$, which satisfies  $ \hat{K}_\alpha (0) = 2^{\alpha -1} \Gamma (\alpha)$ and $(l_t; t \geq 0)$ is a choice of the  local time at level 0 of $R$ made so that the following holds: \\
It\^o's excursion measure $n_\alpha^{\beta \!\downarrow}(de)$ associated with $BES(-\alpha, \beta\!\downarrow)$ together with  this choice of local time may be described as follows:
\begin{itemize}
\item[a)] $ n_\alpha^{\beta \!\downarrow}(V(e)\in dv) = \frac{2^\alpha \Gamma(\alpha+1)}{v^{\alpha +1}} \exp(-\beta v) dv$, where $V(e)$ denotes the lifetime of the generic excursion $e$.
\item[b)] Conditionally on $V=v$, the process $(e(u), u\leq v)$ is distributed as a Bessel bridge of index $\alpha$, and length $v$.
\end{itemize}
Note that this description is valid, in particular, for $\beta = 0$. \\
The conditioned diffusions $BES(-\alpha, \beta\!\downarrow)$, started at $x>0$ and killed at $T_0$ are described in \cite{PY1}. See also Watanabe \cite{W1, W2, W3}.

We now discuss some immediate consequences of these choices: from (a), it follows that:
$$ E_0^{(-\alpha), \beta\!\downarrow} (\exp(-\lambda \tau_u)) = \exp\left(- u 2^{\alpha} \Gamma(\alpha+1)
\int_0^\infty \frac{dv}{v^{\alpha+1}} (1-e^{-\lambda v}) e^{-\beta v} \right)$$ 
and elementary computations yield:
\begin{eqnarray*}
 E_0^{(-\alpha), \beta\!\downarrow} (\exp(-\lambda \tau_u)) &=& \exp\left(- u 2^{\alpha} \Gamma(\alpha+1)\frac{\Gamma(\alpha -1)}{\alpha} \{ (\lambda + \beta)^\alpha - \beta^\alpha\} \right)\\
 &=& \exp\left(- u 2^{\alpha}  \frac{\pi}{\sin(\pi\alpha)} \{ (\lambda + \beta)^\alpha - \beta^\alpha\} \right)
\end{eqnarray*}
As a consequence of these computations, it follows that, on one hand:
\begin{equation} \label{TLstable}
 E_0^{(-\alpha)} (\exp(-\lambda \tau_u)) =  \exp\left(- u 2^{\alpha}  \frac{\pi}{\sin(\pi\alpha)}\lambda ^\alpha  \right)
 \end{equation}
and from \eqref{defbesdrift} taken at $x=0$ and $t= \tau_u$, we can determine the constant $C_\alpha$:
$$C_\alpha = \frac{\pi}{\sin(\pi\alpha)}.$$
On the other hand,
\begin{eqnarray} \label{7}
\lim_{\alpha \vers 0} E_0^{(-\alpha), \beta\!\downarrow} (\exp(-\lambda \tau_u)) &= &
\exp\left(- u  \{ \ln(\lambda + \beta) - \ln(\beta)\} \right)  \nonumber\\
&=& \frac{1}{(1+ \lambda/ \beta)^u}
\end{eqnarray}
Thus, assuming that  we may represent the LHS of \eqref{7} as $E_0^{0, \beta\!\downarrow}(\exp(-\lambda \tau_u))$ for the law $P_0^{0, \beta\!\downarrow}$  of some diffusion $BES(0, \beta\!\downarrow)$, instantaneously reflecting at 0, we will thus have obtained a Krein representation of the Gamma process with parameter $\beta$, i.e. the subordinator whose Laplace transform in $\lambda$, at time $u$, is given by \eqref{7}.
\begin{theo} \label{theo1}
1) For every $\beta >0$, there exists a diffusion $BES(0, \beta\!\downarrow)$ on $\R_+$, which is instantaneously reflecting at 0, and whose infinitesimal generator on $(0,\infty)$ is given by:
$$\frac{1}{2} \frac{d^2}{dx^2} + \left(\frac{1}{2x} + \sqrt{2\beta} \frac{K'_0}{K_0}(\sqrt{2\beta}x) \right) \frac{d}{dx}.$$
2) For some choice of the local time at 0 of $BES(0, \beta\!\downarrow)$, the corresponding It\^o measure of excursions may be described as:
\begin{itemize}
\item[a)] $ n_0^{\beta \!\downarrow}(V(e) \in dv) = \frac{1}{v} \exp(-\beta v) dv$
\item[b)] Conditionally on $V=v$, the process $(e(u), u\leq v)$ is distributed as a Bessel bridge of dimension $2$, and length $v$.
\end{itemize}
3) Let $x >0$, the following relation holds:
 \begin{equation} \label{abscont}
P_x^{(0), \beta\!\downarrow}|_{{\cal R}_t \cap (t<T_0)} = \frac{K_0(\theta R_t)}{K_0 (\theta x)} \exp(  - \beta t) P_x^{(0)}|_{{\cal R}_t}.
\end{equation}
\end{theo}
{\bf Note:} The following relation, where we have extended $(1/K_0)$ to $[0, \infty[$, with $(1/K_0)(0) = 0$,  is equivalent to \eqref{abscont}:
\begin{equation} \label{abscont1}
P_x^{(0)}|_{{\cal R}_t } = \frac{K_0 (\theta x)}{K_0(\theta R_{t\wedge T_0})}\exp(  \beta t) P_x^{(0),  \beta\!\downarrow}|_{{\cal R}_t}.
\end{equation}
The following theorem exhibits an absolute continuity relationship between the laws $P_x^{(0), (\beta + \gamma)\!\downarrow}$ and $P_x^{(0), \beta\!\downarrow}$ from which several important Laplace transforms may be immediately obtained.
\begin{theo} \label{theo2.2}
Let $\beta, \gamma >0$, and $x >0$; then, there is the relationship:
 \begin{equation} \label{abscont2}
P_x^{(0), (\beta+ \gamma)\!\downarrow}|_{{\cal R}_t } = \frac{K_0(\sqrt{2(\beta+\gamma)} R_t)}{K_0 (\sqrt{2\beta} R_t)} \frac{K_0(\sqrt{2\beta} x)}{K_0 (\sqrt{2(\beta+ \gamma)} x)} (1+ \frac{\gamma}{\beta})^{l_t} \exp(  - \gamma t) P_x^{(0),  \beta\!\downarrow}|_{{\cal R}_t}.
\end{equation}
In particular, for $x = 0$ and $t$ being replaced by $ \tau_l$,  formula \eqref{abscont2} simplifies as:
 \begin{equation} \label{abscont3}
P_0^{(0), (\beta+ \gamma)\!\downarrow}|_{{\cal R}_t } = \frac{K_0(\sqrt{2(\beta+\gamma)} R_t)}{K_0 (\sqrt{2\beta} R_t)} (1+ \frac{\gamma}{\beta})^{l_t} \exp(  - \gamma t) P_0^{(0),  \beta\!\downarrow}|_{{\cal R}_t}.
\end{equation}
and
 \begin{equation} \label{abscont4}
P_0^{(0), (\beta+ \gamma)\!\downarrow}|_{{\cal R}_{\tau_l} } =  (1+ \frac{\gamma}{\beta})^{l} \exp(  - \gamma \tau_l) P_x^{(0),  \beta\!\downarrow}|_{{\cal R}_{\tau_l}}.
\end{equation}
\end{theo}
Note, also from \eqref{abscont2}, that the measure:
$$\frac{K_0(\sqrt{2\beta } x)}{K_0 (\sqrt{2\beta} R_t)} \beta^{-l_t} \exp(\beta t) P_x^{(0),  \beta\!\downarrow}|_{{\cal R}_t} $$
does not depend on $\beta$.

\vspace{.3cm}
\noindent
Here are some consequences of \eqref{abscont2},  \eqref{abscont3},  \eqref{abscont4}:
\begin{itemize}
\item[-] As a consequence of  \eqref{abscont2}, one obtains the Laplace transform of \\
$T_0 = \inf \{ t; R_t = 0 \}$ under $P_x^{(0),  \beta\!\downarrow}$:
\begin{equation}
E_x^{(0),  \beta\!\downarrow} (\exp(- \gamma T_0)) =  \frac{K_0 (\sqrt{2(\beta+ \gamma)} x)} {K_0(\sqrt{2\beta} x)}
\end{equation}
In fact, $T_0$ is distributed\footnote{For many references to, and applications of, GIG variables, see Matsumoto-Yor \cite{MY}}  as a $GIG(0; x, \sqrt{2\beta})$ variable, i.e. with density:
$$ \frac{1}{2K_0(x\sqrt{2\beta})} \exp(-\frac{1}{2}(x^2/t + 2 \beta t)) dt.$$ 
\end{itemize}
As a consequence of \eqref{abscont4}, we obtain our main result:
\begin{cor}
$(\tau_l; l \geq 0)$ is, under $P_0^{(0),  \beta\!\downarrow}$, a gamma process with parameter 
$\beta$, i.e. $(\beta \tau_l: l \geq 0)$ is a standard Gamma process:
\begin{equation}
E_0^{(0),  \beta\!\downarrow} (\exp(- \gamma \tau_l)) = \frac{1}{(1+\frac{\gamma}{\beta})^l} \mbox{ for } \gamma \geq 0.
\end{equation}
\end{cor}

%%%%%%%%%%%%%%%%%%%%%%%
%%%%%%%%%%%%%%%%%%%%%%%%%%
\section{Esscher and Girsanov transforms}
This section is devoted to the proofs of the assertions contained in Section 2.

\vspace{.3cm}
\noindent
{\bf (3.1)} To begin with, we explain how, once we have obtained a Krein representation for a subordinator $(S_t)$, we can obtain a related one for a second subordinator $(\tilde{S}_t)$ whose law is an Esscher transform of that of $(S_t)$.
We consider a positive diffusion $(X_t; t \geq 0)$ with 0 as an instantaneous reflecting boundary.
We denote by $L_t$  a choice of its local time at 0 and $(\tau_l: l\geq 0)$ the corresponding inverse local time.
We are looking for a diffusion $\tilde{X}$ on $\R_+$ with laws $\tilde{P}_x$ such that:
\begin{equation} \label{esscher}
\tilde{P}_0|_{{\cal F}_{\tau_u}} = \exp( \psi(a) u - a \tau_u)P_0|_{{\cal F}_{\tau_u}}
\end{equation}
 for $a>0$, where $\psi$ denotes the Laplace exponent of the subordinator $(\tau_u; u \geq0)$.
From \eqref{esscher}, the inverse local time $\tilde{\tau}$ of $\tilde{X}$ is the Esscher transform $\tau^{(a)}$ of $\tau$ defined by:
\begin{equation} \label{Esc}
 E_0[ \exp(- \lambda \tau^{(a)}_u)] = \frac{ E_0[ \exp(- (\lambda +a) \tau_u)] }{E_0[ \exp(-a \tau_u)] }.
 \end{equation}
The  L\' evy measure $m^{(a)}(dx)$ of $\tau^{(a)}$ is related to the L\'evy measure $m(dx)$ of $\tau$ by $m^{(a)} (dx) = \exp(-ax) m(dx)$ (see \cite{Es}, \cite{Sh}, Chapter VII, Section 3c). 

\noindent
To define $\tilde{P}_x$, we need to compute the martingale
$$M_s = E[\exp(-a \tau_t)|{\cal F}_s], \ s \leq \tau_t$$
$$
M_s= \exp(-as)
 E_{X_s(\omega)}\left[\exp( -a\tau_{t-L_s(\omega)})\right]
$$
Now,  we have:
$$
 E_x\left[\exp( -a\tau_v)Ê\right] =
E_x\left[\exp( -aT_0)Ê\right]  E_0\left[\exp( -a\tau_v)Ê\right] $$
and therefore:
$$ \frac{M_s}{M_0} = \varphi_{a\!\downarrow} (X_s) \exp(\psi(a) L_s - as )$$
where $\varphi_{a\!\downarrow} $ denotes the function defined by
\begin{equation}
\varphi_{a\!\downarrow} (x) = E_x[\exp(-a T_0(X))] \qquad (\varphi_{a\!\downarrow} (0) = 1)
\end{equation}
where $T_0(X)$ is the first hitting time of 0 by $X$.
More generally, for $x \geq 0$, we define the law $\tilde{P}_x$ via the absolute continuity relation:
\begin{equation} \label{Girsanov}
\tilde{P}_x|_{{\cal F}_t} = \frac{\varphi_{a\!\downarrow} (X_t)}{\varphi_{a\!\downarrow} (x)} \exp( \psi(a) L_t - a t)P_x|_{{\cal F}_t}.
\end{equation}
The generators $L$  and $\tilde{L}$ of  $X$ and $\tilde{X}$ respectively are linked by:
$$\tilde{L} = L + \frac{\varphi'_{a\!\downarrow} (x) }{\varphi_{a\!\downarrow} (x) } \frac{d}{dx}.$$

\noindent
{\bf (3.2) Examples:} \\
\begin{enumerate}
\item $P_x= P^{(-\alpha)}_x$, $0<\alpha <1$. Then, 
$$\varphi_{a\!\downarrow} (x) = c_\alpha \hat{K}_\alpha(\sqrt{2a}x), \quad \psi(a) =\frac{\pi}{\sin(\pi \alpha)} (2a)^\alpha$$
 and 
$\tilde{P}_x = P_x^{(-\alpha), a\!\downarrow}$. 
\item  $P_x =  P_x^{(-\alpha), \beta\!\downarrow}$, $0<\alpha<1$, $\beta >0$. Then,
$$\varphi_{a\!\downarrow} (x) = c'_\alpha \frac{\hat{K}_\alpha
 (\sqrt{2(\beta+a)}x)}{\hat{K}_\alpha
 (\sqrt{2\beta}x)}; \quad \psi(a) = \frac{\pi}{\sin(\pi \alpha)} \{(2a+ 2 \beta)^\alpha - (2\beta)^\alpha\}$$
 and $\tilde{P}_x = P_x^{(-\alpha), (\beta +a)\!\downarrow}$. 
 \end{enumerate}
 
 \noindent
 From \cite{BY} and \cite{MO}, we know that the inverse local time of $BES(-\alpha)$ is a stable subordinator of index $\alpha$.  It follows from subsection {\bf (3.1)} and Example 1 above that the inverse local time of $BES(-\alpha, \beta\!\downarrow)$ is a subordinator $S^{\alpha, \beta}$. The description of It\^o's excursion measure of $BES(-\alpha, \beta\!\downarrow)$ follows from the description of the excursion measure for $BES(-\alpha)$ given in \cite{BY} and the Esscher transform. Note that  we have not chosen the same normalisation for the local time as in \cite{BY}; we have the following relation  between the two local time processes (from \eqref{TLstable} and subsection 3.2 in \cite{BY}):
 $$l_t^{BY} = \frac{ (2^{\alpha } \Gamma(\alpha))^2}{2(1-\alpha)} l_t ^{DY}$$
 where $l_t^{BY}$ denotes the local time considered in \cite{BY}. Note that there is a mistake in \cite{BY} after (3.i) due to the identification of $\tau^{(-1/2)}_t$ (the inverse local time of the reflected Brownian motion) with $\tau_t$ (for the Brownian motion)  instead of $\tau_{t/2}$  and the correct formula (see p.45, after (3.i)) is 
 $$ E[\exp(-k\tau^{(-\nu)}_{tc(1/2, 2 \nu)})] = \exp(-t c_\nu 2^{\nu -1} k^\nu).$$
 
 \vspace{.3cm}
 \noindent
 {\bf (3.3) Proof of Theorems \ref{theo1} and \ref{theo2.2} } \\
 1) We consider the diffusion $X$ on $[0, \infty[$  with generator
 $$\frac{1}{2} \frac{d^2}{dx^2} + \left(\frac{1}{2x} + \sqrt{2\beta} \frac{K'_0}{K_0}(\sqrt{2\beta}x) \right) \frac{d}{dx}.$$
 A pair $(s(x), m(dx))$ of scale function and speed measure is given by:
 $$s(x) = \int_0^x \frac{dy}{ y K_0^2(\sqrt{2\beta}y)}, \; m(dx) = 2x K_0^2(\sqrt{2\beta}x) dx .$$
 Note that $K_0(y) \sim_0 \ln(2/y)$, thus $\displaystyle \frac{1}{y K_0^2(\sqrt{2\beta}y)}$ is integrable in 0 and $s(0) = 0$. In order to obtain the behavior of $X$ at the boundary point 0, we shall apply the criterium of Rogers-Williams \cite[V.51]{RW} to the diffusion  in natural scale $Y= s(X)$. The speed measure of $Y$ is given by:
 $$m_Y(dx) = \frac{2}{(s'\circ s^{-1}(x))^2} dx.$$
 Then, $\int_{0+} x m_Y(dx) < \infty$ implying that $T_0(Y) < \infty$ $P_x $ p.s. for all $x >0$, and 
$ \int_{0+} m(dx) < \infty$ ensuring that 0 is reflecting. \\
The description of the It\^o measure follows from that of $BES(-\alpha, \beta\!\downarrow)$ letting $\alpha \vers 0$.

\noindent
The absolute continuity relation \eqref{abscont} is given in Pitman-Yor \cite{PY1}. We can also obtain it from \eqref{defbesdrift} and:
$$P_x^{(+\alpha)}|_{{\cal R}_t} = \left(\frac{R_{t\wedge T_0}}{x} \right)^{2\alpha} P_x^{(-\alpha)}|_{{\cal R}_t}, \; x >0.$$
Thus,
$$ \left(\frac{R_{t\wedge T_0}}{x} \right)^{2\alpha} P_x^{(-\alpha), \beta\!\downarrow}|_{{\cal R}_t} 
 = \frac{\hat{K}_\alpha (\theta R_t)}{\hat{K}_\alpha (\theta x)} \exp(  - \beta t) P_x^{(+\alpha)}|_{{\cal R}_t}$$
 and letting $\alpha \vers 0$,
 $$1_{(t < T_0)} . P_x^{(0), \beta\!\downarrow}|_{{\cal R}_t} 
 = \frac{K_0 (\theta R_t)}{K_0 (\theta x)} \exp(  - \beta t) P_x^{(0)}|_{{\cal R}_t}. \qquad \Box $$
 The absolute continuity relation \eqref{abscont2} follows from the example {\bf (3.2)} 2. with $\alpha = 0$.
 In this case,
 $$ \varphi_{a\!\downarrow} (x) = C \frac{K_0(\sqrt{2(\beta+a)}x)}{K_0(\sqrt{2\beta}x)};
  \quad \psi(a) = \ln(1 + \frac{a}{\beta}).$$
The other formulas follow easily. $ \Box$
%%%%%%%%%%%%%%%%%%%%%%%%%%%%%
 \section{Brownian representations of the subordinators $S^{\alpha, \beta}$}
 We keep the same normalisation as in Section 2, i.e. the Laplace transform of a stable process is given by \eqref{TLstable}
\begin{prop}
\begin{enumerate}
\item The stable subordinator $S^{\alpha}$ can be represented by the Brownian additive functional (see
\cite{BY})
$$ S^\alpha_t = A_\alpha (\tau_{\frac{\Gamma(\alpha)^2t}{\alpha^{2\alpha-1}}}) $$
 where
\begin{equation} \label{stabledef}
A_\alpha (s) := \int_0^s |B_u|^{\frac{1}{\alpha}-2} \ du, \; s \geq 0,
\end{equation}
and $\tau_t$ is the inverse local time of $B$.
\item The subordinator $S^{\alpha, 1}$ can be represented as:
\begin{equation} \label{rep}
 \int_0^{\tau_t} h_\alpha  (|B_r|)\  dr.
\end{equation}
where $h_\alpha (x) = 2 (s_\alpha^{-1}(2x))^2 K_\alpha^4(s_\alpha^{-1}(2x))$
and 
\begin{equation} \label{resSL}
s_\alpha(x) = \int_0^x \frac{dy}{y K^2_\alpha(y)}
\end{equation} 
\item The Gamma process $S^{0, 1}$ has the following representation:
\begin{equation}
\left( \int_0^{\tau_t} h(|\beta_u|)\ du\, ; \ t \geq 0 \right)
\end{equation}
with 
\begin{equation} 
 h(x) =  2 (G^{-1} (2x))^2 K_0^4(G^{-1} (2x)),
 \end{equation}
 where $G^{-1}$ is the inverse of the  increasing function $G$ given by:
 \begin{equation} \label{defG}
  G(z) = \int_0^z \frac{dx}{x K_0^2 (x)} \end{equation}
  \end{enumerate}
  \end{prop}
 {\bf Sketch of Proof:} \\ 
 1)The first point  has been obtained in \cite{BY}.
 More precisely,
\begin{equation} \label{stablelaplace}
E[ \exp(- \frac{\lambda}{2} A_\alpha (\tau_t))]  = \exp(-t c_\alpha \lambda^\alpha), \ \lambda \geq 0,
\end{equation}
where 
$$ c_\alpha = \frac{\pi}{\alpha \sin(\pi \alpha)} \left(\frac{\alpha^\alpha}{\Gamma(\alpha)} \right)^2.$$

\noindent
2)  By a Girsanov transform, we can find a diffusion $X$ such that the subordinator
$$A^X_\alpha (\tau_t) := \int_0^{\tau_t(X)} |X_u|^{\frac{1}{\alpha}-2} \ du $$
is a Esscher transform of the stable subordinator $A_\alpha (\tau_t)$ and thus is distributed as $S^{\alpha, 1}$(up to a constant).
Then, we write the diffusion in natural scale $s_\alpha(X)$ as a time change Brownian motion.
We can also start from $BES(\alpha, 1\downarrow)$ and use a time change method.

\noindent
3) By an application of It\^o's formula, we can prove that:
$$ G^{-1} (2 |B_t|) = Y \left( \int_0^t 2 h(|B_s|) ds \right)$$
where $Y$ is a $BES(0, 1\!\downarrow)$ process from which we can deduce that
$\tau_t(Y) \stackrel{(law)}{=} \int_0^{\tau_t (B)} 2 h(|B_s|) ds$
where $\tau_t(Y)$ is the inverse diffusion local time of $Y$, defined with the speed measure
$m_Y(dy) = 2yK_0^2(y) dy$.
   
 \section{Some complements}
 \subsection{The case $\alpha = 1/2$}
 Let us go back to the representation \eqref{rep} for the subordinator with L\'evy measure $k'_\alpha
 \exp(-\frac{y}{2}) \frac{dy}{y^{1+\alpha}}$. In the case $\alpha = 1/2$, \eqref{rep} becomes:
 \begin{equation} \label{rep1/2}
 %A^{1/2}_{1/2}(\tau_t) \stackrel{(law)}{=}
 \int_0^{\tau_t} \frac{dr}{(1+2|\beta_r|)^2}.
 \end{equation}
 We now explain how this representation may be reduced to the more classical one, as given in \eqref{5.2} below.
 For this purpose, we recall that  the stable subordinator of index $1/2$ can also be realized as the distribution of $(T_t; t \geq 0)$ where $T_t$ denotes the first hitting time of $t$ by a Brownian motion, starting from 0. We have the same interpretation for its Esscher transform. Indeed,
 \begin{eqnarray*}
 \ln(1+2|B_t|)& = & 2 \left[ \int_0^t \frac{sgn(B_s) dB_s}{(1+2|B_s|)} -  \int_0^t \frac{ds}{(1+2|B_s|)^2} + L_t\right] \\
 &=&  2 \left[ - \beta^{(1)} \left( \int_0^t \frac{ds}{(1+2|B_s|)^2}\right) + L_t \right]
 \end{eqnarray*}
 where $( \beta^{(1)}(u); u \geq 0) $ is a Brownian motion with drift 1. \\
 From Skorokhod's lemma, we derive:
 $$L_t = \sup\left\{ \beta^{(1)}_s; s \leq  \int_0^t \frac{dr}{(1+2|B_r|)^2} \right\}$$
 Therefore, we obtain:
 \begin{equation} \label{5.2}
 \int_0^{\tau_l} \frac{dr}{(1+2|B_r|)^2} =  \inf \{ u; \beta^{(1)}_u = l \} := T^{(1)}_l.
\end{equation}
Thus, from \eqref{rep1/2}, $(T^{(1)}_l, l \geq0)$ is a subordinator with L\'evy measure $C
 \exp(-\frac{y}{2}) \frac{dy}{y^{3/2}}$ and Laplace transform
 $$ E \left( \exp(- \frac{\lambda^2}{2} T^{(1)}_l) \right) = \exp(-l(\sqrt{\lambda^2 +1}-1)).$$
 Of course, the above Laplace transform could also have  been  computed from the Laplace transform of $T_l$ and the Cameron-Martin formula.
 
 \subsection{Symmetric L\'evy processes}
 From the Brownian representation of the Gamma process, we can give a Brownian representation for the symmetric Gamma process (or variance gamma process) distributed as $(\gamma_1 (t) - \gamma_2(t); t \geq 0)$ where $\gamma_1$ and $\gamma_2$ are two independent gamma processes.
 More generally, we have:
 \begin{prop}
 Let $S_1$ and $S_2$ be two independent subordinators with Brownian representations:
 $$(S_i(t); t \geq 0)  \stackrel{(law)}{=}\left( \int_0^{\tau_t}  ds \varphi_i(|B_s|); t \geq 0\right), i = 1,2.$$
 Then, we have the Brownian representation for $(S_1(t) - S_2(t); t \geq 0)$ as follows:
 \begin{equation} \label{symLev}
 (S_1 (t) - S_2(t); t \geq 0)  \stackrel{(law)}{=} \left( \int_0^{\tau_{2t}}  ds \varphi(B_s); t \geq 0\right)
 \end{equation}
 where $ \varphi (x) = \left\{ \begin{array}{l}  \varphi_1 (x), x >0 \\
 - \varphi_2(- x), x <0 \end{array} \right.$. 
 
 \end{prop}
{\bf Proof:}
It suffices to write
$$ \int_0^{\tau_t}  ds \varphi(B_s) = \int_0^{\tau_t}  ds 1_{(B_s >0)} \varphi_1(B_s) - \int_0^{\tau_t}  ds 1_{(B_s < 0)}\varphi_2(- B_s).$$
Then, we use the representation 
\begin{equation} \label{sym}
 B^+_s = |\beta^{(+)}|_{\int_0^s 1_{(B_u >0)} du} , \;  B^-_s = |\beta^{(-)}|_{\int_0^s 1_{(B_u <0)} du} \end{equation}
where $\beta^{(+)}$ and $\beta^{(-)}$ are (as a consequence of Knight' s theorem) two independent Brownian motions\footnote{ However,  we emphasize that, knowing $B$, only the reflected Brownian motions $|\beta^{(+)}_h|$ and  $|\beta^{(-)}_h|$ are accessible, and not  $\beta^{(+)}$ and  $\beta^{(-)}$.}
; see, e.g., proofs of the arc sine law for Brownian motion, inspired from D. Williams \cite{Wi} (see \cite{KS}, \cite{McK}).\\
Then, we can write 
$$ \int_0^{\tau_t}  ds \varphi(B_s) =  \int_0^{A^{(+)}_{\tau_t}}  dh \varphi_1(|\beta^{(+)}_h| )-
 \int_0^{A^{(-)}_{\tau_t}}  dh \varphi_2(|\beta^{(-)}_h| )$$
 where 
 $$A^{(+)}_{\tau_t} =  \int_0^{\tau_t}  ds 1_{(B_s >0)}  ds; \qquad A^{(-)}_{\tau_t }=  \int_0^{\tau_t}  ds 1_{(B_s <0)}  ds,$$
 and moreover, from \eqref{sym}, we also learn  that:
 $$A^{(+)}_{\tau_t} = \tau_{t/2}(\beta^{(+)}) ,  \;  A^{(-)}_{\tau_t} = \tau_{t/2}(\beta^{(-)}).$$
 Finally, we have obtained the result. $\Box$
 
\subsection{The It\^o measure of $BES(- \alpha, \beta\!\downarrow)$}
The description of the It\^o measure given in Section 2 relies upon one of the descriptions of the It\^o measure $n_\alpha$ of $BES(-\alpha)$ given in \cite{BY}. There is a second description  of $n$ conditionally to the maximum of the excursion.
\begin{prop}
a) Under $n_{\alpha}^{\beta\!\downarrow}$, the distribution of $M$ satisfies:
\begin{equation} 
n_{\alpha}^{\beta\!\downarrow} (M \geq x) = 2 (2 \beta)^\alpha  \frac{ K_\alpha(\sqrt{2 \beta }x)}{ I_\alpha(\sqrt{2 \beta }x)}
\end{equation}
b) Conditionally on $M=x$, the maximum is attained at a unique time $R$ and the processes $(e_t; 0 \leq t \leq R)$ and $(e_{V-t}, 0 \leq t \leq V-R)$ are two independent $BES(\alpha, \beta\!\uparrow)$ processes, starting from 0, stopped at their first hitting time $T_x$ of $x$.
\end{prop}
{\bf Proof:} a) From the description of $n_{\alpha}^{ \beta\!\downarrow}$ given in section 2, we have:
$$ n_{\alpha}^{ \beta\!\downarrow}[f(M^2)] = 2^\alpha \Gamma(\alpha +1) \int_0^\infty \frac{\exp(-\beta v)}{v^{\alpha+1}} \Pi^{(\alpha)}_v [f(M^2)]  dv $$
where $\Pi^{(\alpha)}_v$ is the distribution of a Bessel bridge of index $\alpha$ and length $v$.
$$ n_{\alpha}^{\beta\!\downarrow}[f(M^2)] = 2^\alpha \Gamma(\alpha +1) \int_0^\infty \frac{\exp(-\beta v)}{v^{\alpha+1}} \Pi^{(\alpha)}_1 [f(v M^2)] dv$$
Now, from \cite[Theorem 3.1]{PYIto} 
$$  \Pi^{(\alpha)}_1[\phi(r)] = c_\alpha E[ \phi(\tilde{R}) (\tilde{M})^{-2\alpha}] $$
where the process $\tilde{R}$ is defined by
$$ \tilde{R}_t =  (T + \hat{T})^{-1/2} Y_{t (T + \hat{T})}, \; 0 \leq t \leq 1$$
and
$Y$ is the process connecting the paths of two independent  $BES(2+2\alpha)$ processes $R$ on $[0,T]$ (first hitting time of 1) and $\hat{R}$ on $[0, \hat{T}]$ back to back, i.e.
$$ Y_t = \left\{ \begin{array}{ll} R_t & t\leq T \\
\hat{R}_{T + \hat{T} -t} & T\leq t \leq T+ \hat{T}
\end{array} \right. $$
Then,
 $$\tilde{M} := \sup_{s\leq 1} \tilde{R}_s = (T + \hat{T})^{-1/2}$$
and $c_\alpha =  2^\alpha \Gamma(\alpha +1)$. \\
It follows that:
$$\Pi^{(\alpha)}_1 [f(v M^2)]  = c_\alpha E\left[ f\left( \frac{v}{T + \hat{T}}\right) (T + \hat{T})^\alpha\right].$$
Thus,
\begin{eqnarray*}
 n_{\alpha}^{\beta\!\downarrow}[f(M^2)] &= &2^\alpha \Gamma(\alpha +1) c_\alpha E\left[ \int_0^\infty \frac{\exp(-\beta w(T + \hat{T}) )}{w^{\alpha+1}} f(w) dw \right] \\
 &=& 2^\alpha \Gamma(\alpha +1) c_\alpha  \int_0^\infty \frac{ f(w)}{w^{\alpha+1}}E\left[\exp(-\beta w(T + \hat{T}) ) \right] dw \\
 &=& [2^\alpha\Gamma(\alpha +1)]^{-1} c_\alpha  \int_0^\infty \frac{ f(w)}{w^{\alpha+1}}
 \left( \frac{(\sqrt{2 \beta w})^\alpha}{I_\alpha(\sqrt{2 \beta w})} \right)^2 dw\\
 &=& (2\beta)^\alpha 
  \int_0^\infty \frac{ f(w)}{w I_\alpha^2(\sqrt{2 \beta w})} dw\\
  &=& 2 (2\beta)^\alpha 
  \int_0^\infty \frac{ f(y^2)}{y I_\alpha^2(\sqrt{2 \beta }y)} dy
  \end{eqnarray*}
  %%%%%
  Therefore, we obtain
  $$n_{\alpha}^{\beta\!\downarrow} (M \in dy) = 2 (2 \beta)^\alpha 
  \frac{1}{y I_\alpha^2(\sqrt{2 \beta }y)} dy$$
 hence:
  \begin{eqnarray*}
  n_{\alpha}^{\beta\!\downarrow} (M \geq x ) &= & 2(2 \beta)^\alpha 
 \int_x^\infty  \frac{1}{y I_\alpha^2(\sqrt{2 \beta }y)} dy \\
 &=& 2 (2 \beta)^\alpha  \frac{ K_\alpha(\sqrt{2 \beta }x)}{ I_\alpha(\sqrt{2 \beta }x)}. \qquad 
 \end{eqnarray*}
 We also refer to \cite{PY01} for related computations. \\
 As a verification, we can let $\beta \vers 0$ to obtain:
 $$ n_\alpha (M \geq x) = 2^{2\alpha} \Gamma(\alpha) \Gamma(\alpha +1) \frac{1}{x^{2\alpha}}.$$
 This agrees with the description of the distribution of $M$ under $\hat{n}_\alpha$ given in Biane-Yor \cite{BY}:
$$\hat{n}_\alpha (M \geq x) = \frac{1}{x^{2\alpha}}$$
and from our normalisation given in Section 2, 
$$n_\alpha  = 2^{2\alpha} \Gamma(\alpha) \Gamma(\alpha +1) \hat{n}_\alpha. \qquad $$

\noindent
b) We refer to \cite{PY1} for the definition of Bessel processes with drift  $BES(\alpha,\beta\!\!\uparrow)$.
The description of $n_{\alpha}^{\beta\!\downarrow}$, conditionally to $M$, follows from the description of 
$n_{\alpha}$ (see \cite{BY}), the relation $\displaystyle n_{\alpha}^{\beta\!\downarrow} (de) =
\exp(-\beta V(e)) n_{\alpha}(de)$ and 
$$P_0^{(\alpha), \beta\!\uparrow} |_{{\cal F}_{T_x}} = C_\alpha x^{-\alpha} I_\alpha (\theta x) \exp(-\beta T_x) P_0^{(\alpha)} |_{{\cal F}_{T_x}} $$
(see Proposition 3.1 in \cite{PY1}). $\Box$

 %%%%%%%%%%%%%%%%%%%%%%%%%%%%%%%

%"The Variance Gamma Process and Option Pricing,"
%forthcoming in the European Finance Review (with P. Carr and E. Chang) (1998). 

\end{document}